\documentclass[reqno]{amsart}

\usepackage{latexsym,amssymb,amsmath, amsfonts}
\usepackage{amsthm}
\usepackage{amscd}
\usepackage[shortlabels]{enumitem}
\usepackage{a4wide}

\numberwithin{equation}{section}
\usepackage[usenames,dvipsnames,svgnames,table]{xcolor}
\usepackage{hyperref}
\hypersetup{colorlinks=true,
	linkcolor=blue,
	urlcolor=blue,
	citecolor=red} 
\newcommand{\R}{\mathbb{R}}

\newcommand{\Z}{\mathbb{Z}}

\renewcommand{\le}{\leqslant}
\renewcommand{\ge}{\geqslant}
\renewcommand{\leq}{\leqslant}
\renewcommand{\geq}{\geqslant}

\newcommand{\be}{\begin{equation}}
\newcommand{\en}{\end{equation}}
\newcommand{\ee}{\end{equation}}

\newcommand{\bt}{\begin{theorem}}
\newcommand{\et}{\end{theorem}}
\newcommand{\bp}{\begin{proof}}
\newcommand{\ep}{\end{proof}}
\newcommand{\bc}{\begin{cor}}
\newcommand{\ec}{\end{cor}}
\newcommand{\bl}{\begin{lemma}}
\newcommand{\el}{\end{lemma}}
\newcommand{\bprop}{\begin{prop}}
\newcommand{\eprop}{\end{prop}}

\newcommand{\N}{\mathbb{N}}

\newtheorem{innercustomthm}{Theorem}

\newenvironment{customthm}[1]
{\renewcommand\theinnercustomthm{#1}\innercustomthm}
{\endinnercustomthm}
\newtheorem{theorem}{Theorem}[section]
\newtheorem{remark}{Remark}
\newtheorem{lemma}[theorem]{Lemma}

\newtheorem{proposition}[theorem]{Proposition}

\newtheorem{corollary}[theorem]{Corollary}
\numberwithin{theorem}{section} \numberwithin{definition}{section}

\newcommand{\RNum}[1]{\uppercase\expandafter{\romannumeral #1\relax}}
\def\R{\mathbb{R}}
\def\N{\mathbb{N}}

\def\Z{\mathbb{Z}}

\newcommand{\Rn}{\mathbb{R}^n}

\newcommand{\vertiii}[1]{{\left\vert\kern-0.25ex\left\vert\kern-0.25ex\left\vert #1 
		\right\vert\kern-0.25ex\right\vert\kern-0.25ex\right\vert}}

\theoremstyle{definition}
\usepackage[leqno]{amsmath}
\makeatletter
\newcommand{\leqnomode}{\tagsleft@true}
\newcommand{\reqnomode}{\tagsleft@false}
\makeatother

\author[A. Mu\~noz]{Alexander Mu\~noz}
\address{IME - Universidade de São Paulo, R. do Matão, 1010 - Butantã, São Paulo - SP, Brazil.}
\email{alexd@ime.usp.br}

\thanks{}

\date{}
\title[Decay of solutions to some KdV perturbations]{On decay of solutions to some perturbations of the Korteweg-de Vries equation}

\begin{document}

\begin{abstract} 
This work is devoted to study the relation between regularity and decay of solutions of some dissipative perturbations of the Korteweg-de Vries equation in weighted and asymmetrically weighted Sobolev spaces.
\end{abstract}

\maketitle

\section{Introduction}\label{intro}
In this work, we explore dissipative perturbations of the Korteweg-de Vries equation posed over weighted spaces, such as $H^s(\R)\cap L^2(|x|^{2r}dx)$ or $H^s(\R)\cap L^2(e^{2rx}dx)$. We adopt the classical notion of well-posedness established by Kato (\cite{KATO}), where an initial value problem (IVP) is considered locally well-posed in a Banach space $Y$ if, for every initial datum $f_0\in Y$, there exist $T>0$ and a unique solution in the space $C([0,T];Y)$ with the data-solution map $f_0\mapsto f$ continuous from $Y$ to $C([0,T];Y)$.

Motivated by the question of whether dispersive equations preserve the flow of the Schwartz space, several classical works have addressed the study of dispersive equations over spaces that measure, in addition to regularity, the decay of solutions. In \cite{KATO}, T. Kato investigated the Korteweg-de Vries (KdV) equation
\begin{equation*}
\partial_t u + \partial_x^3 u + u\partial_x u=0
\end{equation*}
in weighted Sobolev spaces $Z_{s,r}:=H^s(\R)\cap L^2(|x|^{2r}dx)$ with integer indices $s, r$ and asymmetrically weighted Sobolev spaces $H^s(\R)\cap L^2(e^{2rx}dx)$ where $r>0$. When considering $Z_{s,r}$ spaces, numerous works have extended Kato's results to settings where the indices are allowed to be fractional (see for instance \cite{FLP}). These ideas were further generalized to other dispersive models, such as the two-dimensional Zakharov-Kustensov (ZK) equation:
\begin{equation*}
\partial_t u + \partial_{x_1} \Delta u + u\partial_{x_1}u=0
\end{equation*}
The study of this equation involved classical \cite{BJM} and anisotropic weights \cite{pachon}. A similar situation occurred for other well-known dispersive models, like the Schrödinger equation. The seminal works of Hayashi, Nakamitsu, and Tsutsumi (\cite{Tsut1}, \cite{Tsut2} and \cite{Tsut3}) were extended in \cite{NAHASPONCE} to fractional powers.

In recent years, a substantial amount of research has focused on studying the decay of solutions to a wide variety of dispersive models in weighted spaces. See, for instance, \cite{Ostrosky}, \cite{BO}, \cite{MP1} and the references therein.

In the study of asymmetrically weighted Sobolev spaces, Kato's work has recently been extended to the context of the ZK equation in \cite{felipeexp}, where the authors established local well-posedness in dimensions $n\ge 2$ using energy estimates. Additionally, in \cite{ILP2013}, the authors conducted a study measuring how this exponential property decreases with time. \cite{ILP2013} also introduced a notion of optimality regarding the behavior of the indices $s$ and $r$ for weighted Sobolev spaces. Specifically, they established that if the solution in $L^2(\R)$ of the KdV equation exhibits an excess of decay at two different times, then the regularity of the solution improves proportionally to the excess of decay. A similar study was conducted in \cite{5kdv} for the fifth-order KdV equation.

More recently, in \cite{MP2}, these notions of optimality were revisited for the modified KdV equation and improved in the sense that arbitrarily large excess of decay is accepted. This improvement allows this effect to occur in $H^s$ for any $s\ge0$.

This work aims to conduct a comprehensive study of the decay properties of several dispersive equations that also exhibit a dissipative nature. We delve into the local well-posedness in $H^s(\R)$ and $Z_{s,r}$ spaces. Additionally, we study the optimality of the relation between $s$ and $r$, along with the local well-posedness in asymmetrically weighted spaces

Let us specify the exact problem we will address. Let $\eta>0$ and define $\widehat{Lu}(\xi)=-\Phi(\xi)\widehat{u}(\xi)$ where $\Phi(\xi)=-|\xi|^p+\Phi_1(\xi)$ with $p\in \R^+$ and $|\Phi_1(\xi)|\le c(1+|\xi|^q)$ for some $0\le q<p$. We are interested in the initial value problem (IVP)
\begin{equation}\label{P}
    \begin{cases}
        \partial_t u+\partial^3_x u + \eta L u + u\partial_x u = 0 \ \ t>0, \ x\in \R, \\
        u(x,0)=u_0(x).
    \end{cases}
\end{equation} 

Some examples of models that fit into the definition of \eqref{P} are, for instance, the Korteweg-de Vries-Burguers (KdVB) equation $$\partial_t u +\partial_x^3u -\eta \partial_x^2 u+u\partial_x u=0,$$ the Ostrovsky-Stepanyants-Tsimring (OST) equation $$\partial_t u +\partial_x^3 u -\eta \mathcal{H}(\partial_x u+\partial_x^3)+u\partial_x u=0$$ or the Korteweg-de Vries-Kuramoto-Sivashinsky (KdV-KS) equation $$\partial_t u+\partial_x^3 u+\eta(\partial_x^2+\partial_x^4)+u\partial_x u=0.$$ For more details about these equations see, for instance, \cite{ost}, \cite{burgers} and \cite{kdvks}.

The local well-posedness in $H^s(\R)$ of IVP \eqref{P} was established in \cite{critpert} for $s>-p/2$. More precisely, it was established the following result.
\begin{customthm}{A}
    Let $\eta>0$ and $p>3$. For any data $u_0\in H^s(\R)$ with $s>-p/2$ there exists a time $T(\|u_0\|_{H^s})$ and a unique solution $v$ to the IVP \eqref{P} in $C([0,T];H^s(\R))$. Moreover, the map $v_0\mapsto v$ is smooth from $H^s(\R)$ to $C([0,T];H^s(\R))\cap \cdots$.
\end{customthm}
The authors also prove their result is sharp in the sense for $p\ge2$ and $s<-p/2$ the data-solution map cannot be smooth and therefore no solution can be obtained from the contraction mapping principle for initial data with Sobolev regularity below $-p/2$. As far as we know, the study of \eqref{P} in weighted spaces has not been addressed in a general context. To the best of our knowledge, there is only one particular result facing a member of this family under weighted Sobolev spaces and it was done in \cite{MP1} for the OST equation. The authors proved local well-posedness in $Z_{s,r}$ for $s>0$ and $r\le s/2$. 

\subsection{Main Results}\hfill

Our first result is related to the local well-posedness in $H^s(\mathbb{R})$. Despite the estimates provided by Carvajal and Panthee in \cite{critpert} can be modified to allow polynomial decay, the computations involving the time variable seem to hinder such modifications to recover the same existence time $T$ when studying local well-posedness in $Z_{s,r}$ with $s>0$. This limitation prevents us from obtaining a persistence result (as in Theorem \ref{sametime} below), which is crucial when dealing with the optimal relation between $s$ and $r$. Therefore, in Theorem \ref{LWPHS} we introduce an alternative local well-posedness in $H^s(\R)$ for $s\ge 0$. 

Based on the Banach space defined in the proof of Theorem \ref{LWPHS}, we can include extra norms that ensure decay properties to the respective solution provided by the contraction principle. In Theorem \ref{wLWP} we establish local well-posedness $Z_{s,r}$ for $s>0$ and $r\le \frac{s}{p-1}$ in the interesting case $p>3$. Note if $p\le3$, according to the definition of $K$ in Theorem \ref{troca} (below), the decay would be controlled by $K=2$ coming form $\partial_x^3$ instead of $L$. 

Following the ideas of \cite{NAHASPONCEprep}, in Theorem \ref{sametime} we show that the result in $Z_{s,r}$ can be seen as a persistence result in which the solution in $H^s(\R)$ persists in $Z_{s,r}$ when the initial datum also belong to $L^2(|x|^{2r}dx)$. 

Our main result aims to study the optimal relationship between $s$ and $r$ obtained from the well-posedness in weighted Sobolev spaces, in the sense of finding a upper bound for decay $r\le s/\beta$ which one cannot expect to be relaxed. In fact, an excess of decay $s/\beta+\varepsilon$ translates into an improvement of regularity proportional to the excess $s+\beta\varepsilon$. We focus on a specific subfamily of equations considered in the definition of $\Phi$, which comprises an even leading power $p=2k$, and $\Phi_1$ as a linear combination of terms in the form $\xi^m|\xi|^n$, where $m$ is odd, $n$ is even, and $n+m<2k$. In this case, it turns out that the relationship $r\le\frac{s}{p-1}=\frac{s}{2k-1}$ is not optimal as it was expected from Theorem \ref{troca} and the works available in the literature (see \cite[Theorem 1.6]{5kdv}, \cite[Theorem 1.4]{ILP2013} and \cite[Theorem 1.1]{MP2}). We proved in Theorem \ref{optimal} that the optimal relation is actually $r\le \frac{s}{2k}$ rather than the bound $r\le\frac{s}{2k-1}$ provided by Theorem \ref{wLWP}. The latter shows there is a stronger gain in regularity when the decay exceeds $s/2k$ and shows that for this particular family of equations the regularity in Theorem \ref{wLWP} is actually greater when $s/2k<r\le \frac{s}{2k-1}$ (see Corollary \ref{coro} below).

The final result of this paper, Theorem \ref{exp}, shows the local well-posedness of \eqref{P} in asymmetrically weighted Sobolev spaces $H^s(\R)\cap L^2(e^{2rx}dx)$ for any $r>0$ where $p>1$ is an even integer and $\Phi_1$ taken as a linear combination of terms of the form $\xi^m|\xi|^{2n}$ with $n,m$ integers so that $m+2n<p$ .  \vspace{2mm}

This paper is organized as follows: In Section \ref{prelim}, we introduce notation and develop the preliminary theory related to derivatives, weights, and linear estimates. Section \ref{sec3} presents the proof of local well-posedness results in $H^s$ and $Z_{s,r}$, along with a persistence result in $Z_{s,r}$. In Section \ref{sec4}, we establish the optimality of the relationship between regularity and decay, in the sense discussed earlier. Finally, Section \ref{sec5} provides the proof of a local well-posedness result in asymmetric spaces with exponential weights.

\section{Preliminaries}\label{prelim}
\subsection{Notation}\hfill

We use $c$ to denote several constants that may vary from line to line. Dependence on parameters is indicated using subscripts or parenthesis. 

$\mathcal{S}(\R^n)$ denotes the Schwartz space, $L^p(\R)$ is the usual Lebesgue space and $L^2(wdx)$ represent the space $L^2$ with respect to the measure $w(x)dx$. The norm in $L^p(\R)$ will be denoted by $\|\cdot\|_p$. For a function $f$, the expressions $\widehat{f}$ and $f^\vee$ mean the Fourier and inverse Fourier transform, respectively.  For $s, b\in \R$ the set $H^s(\R)$ is the $L^2$-based Sobolev space of index $s$ and $Z_{s,b}$ is the weighted Sobolev space $H^s(\R)\cap L^2(|x|^{2b}dx)$. 
We denote with $\|\cdot\|_{L^p_xL^q_T}$ the norm in the mixed Lebesgue spaces defined by \begin{equation*}
    \|f\|_{L^p_xL^q_T}:= \left( \int_\mathbb{R} \left( \int_0^T |f(x,t)|^q dt\right)^{p/q} dx
\right)^{1/p}.
\end{equation*}
When the integration in time is over $\R$ we write $\|f\|_{L^p_xL^q_t}$ and in case $p=q$, we have $\|f\|_{L^p_xL^q_T}=\|f\|_{L^q_{T}L^p_x}=:\|f\|_{L^p_{xT}}$. We use $D^s_x$ for the classical fractional derivative $\widehat{D_x^sf}(\xi)=|\xi|^s \widehat{f}(\xi)$. For simplicity we adopt the notation $\langle x\rangle:=(1+x^2)^{1/2}$. Also, denote with $J^sf$ the potential $\widehat{J^s(f)}(\xi)=\langle \xi \rangle^{s}\widehat{f}(\xi).$ Given a real number $x$ we use $x^+$ (respect. $x^-$) to mean $x+\varepsilon$ (respect. $x-\varepsilon$) for some arbitrarily small $\varepsilon>0$. 

\subsection{Preliminary estimates}

We begin by stating several results related to derivatives and weights that will be useful in the proofs from Section \ref{sec3} to Section \ref{sec5}. The following result is a general estimate related to the Kato-Ponce inequality, as established by Grafakos and Oh.

\begin{theorem}\label{graf}
    Let $1/2<r<\infty$, $1<p_1,p_2,q_1,q_2\le \infty$ satisfy $\frac{1}{r}=\frac{1}{p_1}+\frac{1}{q_1}=\frac{1}{p_2}+\frac{1}{q_2}$. Given $s>\max\{0, \frac{n}{r}-n\}$ or $s\in 2\N$, there exists $c=c(n,s,r,p_1,q_1,p_2,q_2)<\infty$ such that for all $f, g \in \mathcal{S}(\R^n)$ we have
    \begin{equation*}
        \|D^s(fg)\|_{L^r}\le c\left\{\|D^sf\|_{L^{p_1}}\|g\|_{L^{q_1}}+\|f\|_{L^{p_2}}\|D^sg\|_{L^{q_2}} \right\}.
    \end{equation*}
\end{theorem}
\begin{proof}
    See Theorem 1 in \cite{grafakosoh}
\end{proof}

Denote with $A_p$ the Muckenhoupt class on $\R^n$. More precisely, given $1<p<\infty$, $A_p$ consists of all weights $\omega$ such that
\begin{equation}
	[\omega]_p=\sup_{Q}\left(\frac{1}{|Q|}\int_{Q}\omega(y)dy \right)\left(\frac{1}{|Q|}\int_{Q}\omega^{-\frac{1}{p-1}}(y)dy \right)^{p-1}<\infty,
\end{equation}
where the supremum is taken over all cubes $Q\subset\R^n$; additional details may be seen in \cite{Hunt}.

\begin{lemma}\label{hilberpp}
The Hilbert transform is bounded in $L^p(wdx)$, $1<p<\infty$, if and only if $w\in A_p$.\vspace{-2mm}
\end{lemma}
\begin{proof} See Theorem 9 in \cite{Hunt}.
\end{proof}

The following lemma deals with the interplay between regularity and decay and will be used often as an "interpolating" step in the upcoming computations.

\begin{lemma}\label{interplemma}
	Assume $a, b>0$, $1<p<\infty$ and $\theta\in(0,1)$. If $J^a f\in L^p(\R^n)$ and $\langle x \rangle^b f\in L^p(\R^n)$ then 
	\begin{equation}\label{interp2}
		\|\langle x\rangle^{\theta b}J^{(1-\theta) a}f\|_{L^p(\R^n)}\le C \|\langle x\rangle^b f\|_{L^p(\R^n)}^{\theta}\|J^a f\|_{L^p(\R^n)}^{1-\theta}.
	\end{equation}
	The same holds for $D$ instead of $J$. Moreover, for $p=2$ we have
	\begin{equation}\label{interp}
		\left\|J^{\theta a}\left(\langle x\rangle^{(1-\theta)b}f\right)\right\|_{L^2(\R^n)}\le C \|J^a f\|_{L^2(\R^n)}^\theta \|\langle x\rangle^b f\|_{L^2(\R^n)}^{1-\theta}.
	\end{equation}
\end{lemma}
\begin{proof}
	Inequality	\eqref{interp} follows from \eqref{interp2} in view of Plancherel's identity. The proof of \eqref{interp2} follows using Hadamard's three lines theorem. See Lemma 4 in \cite{NAHASPONCE}.
\end{proof}
We continue by introducing a result regarding the exchange between the linear group associated to the equation and polynomial weights.

For a continuous real-valued function $\phi$ consider the regularity conditions  
\leqnomode
\begin{equation}\tag{A}
	\begin{split}
		&\text{ There exists a continuous function } g:\R^n\to\R,  \ g>0 \text{ except maybe at }x=0,\\&\text{ so that for all }
		x,y\in \Rn\text{ with }|x-y|\le |x|\text{ we have } |\phi(x)-\phi(y)|\le g(x) |x-y|. \label{AA}
	\end{split}
\end{equation}
\begin{equation}\tag{B}
	\begin{split}
		&\text{There exists } C>0 \text{ such that for all }  x,y\in\Rn \text{ satisfying } |x-y|\ge|x| \text{ we have}\\& |\phi(x)-\phi(y)|\le C |x-y|^a, \text{for some }a\ge 1.  \label{BB}
	\end{split}
\end{equation}
\reqnomode

\begin{theorem}\label{troca}
	Let $m\in \Z^+$ and assume that $\phi_1,\dots, \phi_m$ satisfy conditions \eqref{AA} and \eqref{BB} with $g_i(x)\le c_i (1+|x|^{k_i})$, for some $k_i\in\Z^+$ and $c_i>0$, $i=1\ldots,m$. Set $$\Phi(\xi):=\displaystyle\sum_{i=1}^{m}\phi_i(\xi)$$ and $K:=\max\limits_i k_i$. Let $L$ be the linear operator defined by $Lf=\left(i\Phi(\xi)\widehat{f} \right)^\vee$ and assume $0<s<K$.
	If $u\in C([-T,T],H^s(\R^n))$ is the solution of the IVP
	\begin{equation}\label{lpvi}
		\begin{cases}
			\partial_t u +Lu=0, \quad x\in\R^n, \ t\in \R,  \\
			u(0)=u_0\in Z_{s,r}:=H^s(\R^n)\cap L^2(|x|^{2r}dx),
		\end{cases}
	\end{equation}
	with  $0<r\le s/K$, then $u$ satisfies the inequality
	\begin{equation}\label{trocas}
		\||x|^r u(t) \|_{L^2}=\||x|^r U(t)u_0 \|_{L^2}\le C\left\{(1+|t|)\|u_0\|_{H^s} + \||x|^ru_0\|_{L^2}  \right\},
	\end{equation}
	where $c$ depends on $K$, $m$, $s$ and $n$.
\end{theorem} 
\begin{proof}
    See Theorem 1.1 in \cite{MP1}.
\end{proof}
The restriction $s<K$ in Theorem \ref{troca} can be removed via iterations, see Section 1.5 in \cite{tese} for the details.
\subsection{Linear estimates}
The solution to the linear part of \eqref{P} is given by $V(t)u_0(x)$ where the group $V(t)$ is defined by
\begin{equation}
    V(t)u_0(x):=(e^{it\xi^3+\eta t\Phi(\xi)}\widehat{u}_0(\xi) )^\vee.
\end{equation}

Following the ideas of \cite{CarvScia}, in \cite{critpert}, the authors introduced a refined local smoothing effect for a phase function $\Phi$ of the form $|\xi|^k - |\xi|^{k+2}$, where $k \in \mathbb{Z}^+$. A similar result can be obtained for the phase function defined in \eqref{P}. For the sake of completeness, we provide the details here.

Note from the definition of $\Phi$, there exists $M>0$ large enough such that for all $|\xi|\ge M$ we have
\begin{equation}\label{dif}
    |\Phi(\xi)|\ge \frac{|\xi|^p}{2} \mbox{  and  } \Phi_1(\xi)\le \frac{|\xi|^p}{2}.
\end{equation}
Moreover, it can be seen that for $0<t<T$ we have $\|e^{\eta t\Phi(\xi)}\|_{L^\infty_\xi}\le c(M)e^{\eta T}$.

\begin{theorem}
    \label{smooth}
    Let $T>0$, $u_0\in L^{a_1}(\mathbb{R})$, $0\leq s < \frac{p}{a}-\frac{1}{{a_1}}+\frac{1}{b}$, $a\geq 2$, ${a_1}\geq 2$, $b\geq 2$ satisfying $1=\frac{1}{a}+\frac{1}{{a_1}}$.
    We have
    \begin{equation*}
        \| D_x^s V(t) u_0 \|_{L_T^aL_x^b}\leq c(\eta, M) \{ e^{\eta T}T^{1/a} +(a\alpha)^{-1/a}T^\alpha  \} \|  u_0 \|_{L_x^{a_1}}, 
    \end{equation*}
    where $\alpha= \alpha(a,b,s):=\frac{1}{a}-\frac{s}{p}-\frac{1}{p} \left(  \frac{1}{{a_1}}-\frac{1}{b} \right) $
\end{theorem}

Particular cases of Theorem \ref{smooth} will be used so frequently that we choose to state them separately.

\begin{corollary}
    \label{cor1}
    Let $T>0$, $u_0\in L^{a_1}(\R)$, $0\leq s < \frac{p}{a}-\frac{1}{{a_1}}$, $a\geq 2$, ${a_1}\geq 2$ satisfying $1=\frac{1}{a}+\frac{1}{{a_1}}$.
    We have 
    \begin{equation*}
        \| D_x^s V(t)u_0  \|_{L_T^a L_x^\infty} \leq c(\eta, M) \{ e^{\eta T}T^{1/a}+(a\alpha)^{-1/a}T^{\alpha} \} \| u_0 \|_{L_x^{a_1}},
    \end{equation*}
    where $\alpha= \alpha(a,\infty, s) = \frac{1}{a}-\frac{s}{p}-\frac{1}{p{a_1}}.$
\end{corollary}

\begin{corollary}
    \label{cor2}
    Let $T>0$, $u_0\in L^2(\mathbb{R})$, $b\geq 2$ and $0\leq s < \frac{p-1}{2}+\frac{1}{b}$. We have
    \begin{equation*}
        \| D_x^s V(t)u_0\|_{L_T^2L_x^b} \leq c(\eta, M) \{ e^{\eta T}T^{1/2}+(2\alpha)^{-1/2}T^{\alpha} \} \| u_0 \|_{L_x^2},
    \end{equation*}
    where $\alpha=\alpha(2,b,s)=\frac{1}{2}-\frac{s}{p}-\frac{1}{p} \left( \frac{1}{2} -\frac{1}{b} \right)$.
\end{corollary}

\begin{corollary}
    \label{cor3}
    Let $u_0 \in \dot{H}^s(\mathbb{R})$, $T>0$, $b\geq 2$ and $\frac{3}{2}-\frac{1}{b}-\frac{p}{2}$. We have 
    \begin{equation*}
        \| D_x^1 V(t)u_0 \|_{L_T^2L_x^b} \leq c(\eta, M)\{ e^{\eta T}T^{1/2} +(2\alpha)^{-1/2}T^{\alpha} \} \| D_x^s u_0 \|_{L_x^2}, 
    \end{equation*}
    where $\alpha=\alpha(2,b,1-s)=\frac{1}{2}-\frac{1}{p}+\frac{s}{p}-\frac{1}{p}\left(\frac{1}{2}-\frac{1}{b} \right)$.
\end{corollary}

\begin{corollary}
    \label{cor4}
    Let $T>0$ $u_0\in L^2(\mathbb{R})$ with $0\leq s <p/2$. We have
    \begin{equation*}
        \| D_x^s V(t)u_0 \|_{L_{xT}^2}\leq c(\eta, M) \{ e^{\eta T}T^{1/2}+(2\alpha)^{-1/2}T^{\alpha} \} \| u_0 \|_{L_x^2},
    \end{equation*}
    where $\alpha=\alpha(2,2,s)=\frac{1}{2}-\frac{s}{p}$.
\end{corollary}

\begin{proof}[Proof of Theorem \ref{smooth}]

Let $u_0\in \mathcal{S}(\mathbb{R})$. Consider the cut-off function $\varphi \in C(\mathbb{R} \setminus \{ 0 \})$, $0\leq \varphi \leq 1$, defined by 
\begin{equation*}
        \varphi (t)=
\begin{cases}
    1,& 0\leq t\leq 1,\\
    0,& t\in [0,2)^c.   
\end{cases}
\end{equation*}

Set $\varphi_T(t) := \varphi(t/T).$ We have that 
\begin{align*}
    \| D_x^s V(t)u_0 \|_{L_T^aL_x^b} 
    & = \left( \int_0^T \left( \int | D_x^s V(t) u_0(x)|^bdx \right)^{a/b}  dt \right)^{1/a}\\
    &\leq \left( \int_{-\infty}^\infty \left( \int |\varphi_T(t) D_x^s V(t) u_0(x) |^b dx \right)^{a/b} dt \right)^{1/a}\\
    &=\| \varphi_T(\cdot) D_x^s V(t)u_0 \|_{L_t^aL_x^b}.
\end{align*}

By duality, it is enough to prove that for all $g\in L_t^{a_1}L_x^{b_1}$ with $1=\frac{1}{a}+\frac{1}{{a_1}}=\frac{1}{b}+\frac{1}{{b_1}}$, it follows
\begin{equation*}
    \left| \int\int \varphi_T(t) D_x^s V(t)u_0(x) g(x,t)dxdt \right| \leq c\| u_0 \|_{L_x^{a_1}}\| g\|_{L_t^{a_1}L_x^{b_1}}.
\end{equation*}

Note 
\begin{align*}
    D_x^sV(t)u_0(x)&= \left( |\xi|^se^{it\xi^3+\eta t \Phi(\xi)}\hat{u}_0 \right)^{\vee}(x)\\
    & = c\int e^{ix\xi}|\xi|^s e^{it\xi^3+\eta t\Phi(\xi)}\hat{u}_0(\xi)d\xi,
\end{align*}
thus,
\begin{align}
    \left| \int \int \varphi_T(t)D_x^sV(t)u_0(x)g(x,t)dxdt \right|
    &\leq c\left| \int \int \int \varphi_T(t)e^{ix\xi +it\xi^3+\eta t\Phi(\xi)}|\xi|^s \hat{u}_0(\xi) g(x,t) d\xi dxdt \right| \nonumber\\
    &\leq c \| \hat{u}_0 \|_{L_\xi^a} \left\| |\xi|^s \int \varphi_T(t) e^{\eta t \Phi(\xi)} \int e^{ix\xi}g(x,t)dxdt \right\|_{L_\xi^{a_1}}\nonumber\\
    &\leq c\| u_0 \|_{L_x^{a_1}} \left\| |\xi|^s \int \varphi_T(t)e^{\eta t \Phi(\xi)}g^\vee(\cdot,t)(\xi)dt \right\|_{L_\xi^{a_1}}\nonumber\\
    &=: c\| u_0 \|_{L_x^{a_1}} \| Ag \|_{L_\xi^{a_1}}.\label{3b}
\end{align}

First, note that for $\gamma$ such that $\frac{1}{{a_1}}=\frac{1}{\gamma}+\frac{1}{b}$ we have
\begin{align}
    \| Ag \|_{L_\xi^{a_1}(|\xi|\leq M)} 
    &\leq c(M) \int \varphi_T(t) \left\| e^{\eta t\Phi(\xi)}g^\vee(\cdot,t) \right\|_{L_\xi^{a_1}(|\xi|\leq M )}dt\nonumber\\
    & \leq c(M) \int \varphi_T(t) \| e^{\eta t \Phi(\xi)} \|_{L_\xi^\gamma(|\xi|\leq M)} \| g^\vee (\cdot,t) \|_{L_\xi ^b (|\xi|\leq M)}dt \nonumber \\
    &\leq c(M) e^{\eta T}\| \varphi_T \|_{L_t^a} \left\| \|g(\cdot,t)\|_{L_x^{b_1}} \right\|_{L_t^{a_1}}\nonumber\\
    &\leq c(M)e^{\eta T}T^{1/a}\| g \|_{L_t^{a_1}L_x^{b_1}}.\label{3}
\end{align}

On the other hand, 
\begin{align*}
    \| Ag \|_{L_\xi^{a_1} (|\xi|\geq M)} &
    \leq \int \varphi_T(t) \left\| |\xi|^se^{\eta t \Phi(\xi)} \right\|_{L_\xi^\gamma (|\xi|\geq M)} \| g(\cdot, t) \|_{L_x^{b_1}}dt.
\end{align*}

According to \eqref{dif},
\begin{align*}
    \left\| |\xi|^se^{\eta t\Phi(\xi)} \right\|_{L_\xi^\gamma(|\xi|\geq M)}
    &\leq \left\| |\xi|^s e^{-\frac{\eta t}{2}|\xi|^p} \right\|_{L_\xi^\gamma(|\xi|\geq M)}
    \leq c(\eta,M)t^{-\left(\frac{s}{p}+\frac{1}{pr}\right)}.
\end{align*}

Hence
\begin{align}
    \| Ag \|_{L_\xi^{a_1}(|\xi|\geq M)} 
    &\leq c(\eta,M) \left\| \varphi_T(t)t^{-\left( 
s/p+1/pr \right)} \right\|_{L_t^a}\| g \|_{L_t^{a_1}L_x^{b_1}} \nonumber\\
&\leq c(\eta,M)(a\alpha)^{-1/a}T^\alpha\| g \|_{L_t^{a_1}L_x^{b_1}}\label{4},
\end{align}
where $\alpha = \frac{1}{a}-\left( \frac{s}{p}+\frac{1}{pr} \right).$

From \eqref{3} and \eqref{4} we conclude 
\begin{align}
\label{5}
    \| Ag \|_{L_\xi^{a_1}}\leq c(\eta, M) \{ e^{\eta T}T^{1/a}+(a\alpha)^{-1/a}T^\alpha \} \| g \|_{L_t^{a_1}L_x^{b_1}}.
\end{align}

The result follows from \eqref{3b}, \eqref{3} and \eqref{5}.
\end{proof}

\begin{proposition}
    \label{infinity}
    For $T,p,q>0$ real numbers with $p>2q$. We have 
    \begin{equation}
        \label{infinityeq}
        \|D^q_xV(t)u_0\|_{L^\infty_xL^2_T}\le c(\eta, M, T, p, q)\|u_0\|_{L^2_x}. 
    \end{equation}
\end{proposition}
\begin{proof}
    Note for $x\in \R$
    \begin{equation*}
    \begin{split}
        D^q_xV(t)u_0(x)=c\int e^{ix\xi}e^{it\xi^3+\eta t\Phi(\xi)}|\xi|^q\widehat{u}_0(\xi)d\xi.
    \end{split}     
    \end{equation*}
Therefore, according to Plancherel's identity
\begin{equation*}
    \begin{split}
        \|D^q_xV(t)u_0(x)\|_{L^2_T}^2&=\int_0^T\left| \int e^{ix\xi}e^{it\xi^3+\eta t\Phi(\xi)}|\xi|^q\widehat{u}_0(\xi)d\xi  \right|^2dt\\&\le \int_0^T\int e^{2\eta t\Phi(\xi)}|\xi|^{2q}|\widehat{u}_0(\xi)|^2d\xi dt\\&\le \int_0^T\|e^{2\eta t\Phi(\xi)}|\xi|^{2q}\|_{L^\infty_\xi}\|u_0\|^2_{L^2_x}dt.
    \end{split}
\end{equation*}
Using to \eqref{dif}, it can be seen that $\|e^{2\eta t\Phi(\xi)}|\xi|^{2q}\|_{L^\infty_\xi}\le c(\eta,M,p,q)(1+t^{-\frac{2q}{p}})$. Thus,
\begin{equation}
    \|D^q_xV(t)u_0(x)\|_{L^2_T}^2\le c(\eta,M,p,q)\|(1+t^{-\frac{2q}{p}})\|_{L^1_T}\|u_0\|_{L^2_x}\le c(\eta,M,p,q,T)\|u_0\|_{L^2_x}.
\end{equation}

\end{proof}

\subsection{Integration formulae}

We introduce formal computations that allow us to obtain integration patterns useful for performing energy-type estimates.

We define, for $u$ and $\rho$ sufficiently well-behaved to perform integration by parts, the bracket
\begin{align*}
    \langle n,m,a \rangle := \int \partial_x^n u\partial_x^m u \rho^{(a)}dx.
\end{align*}
\begin{proposition}
    Let $n\in\mathbb{N}$. For any $k\in \{ 1,2,\dots ,n \}$ we have
    \begin{align*}
        \langle n,n-k,a \rangle = 
        \sum_{j=0}^{\lfloor k/2 \rfloor} c_j(-1)^{k+j} \langle n-k+j, n-k+j, a+k-2j \rangle,
    \end{align*}
    where $a\in\mathbb{Z}^+ \cup \{ 0 \}$, $c_j>0$ and $\lfloor x \rfloor$ stands for the greatest integer less than or equal to $x$.
\end{proposition}

\begin{proof}
    The proof is done by induction over $n$.

    \textbf{Base case: $n=2$}.
    
    \noindent If $k=1$, integrating by parts
    \begin{equation*}
        \langle 2,1,a \rangle = -\frac{1}{2} \langle 1,1,a+1 \rangle.
    \end{equation*}
    If $k=2$, integrating by parts
    \begin{equation*}
        \langle 2,0,a \rangle = -\langle 1,1,a \rangle - \langle 1,0,a+1 \rangle,
    \end{equation*}
    where $-\langle 1,0,a+1 \rangle = \frac{1}{2} \langle 0,0,a+2 \rangle$.
Thus 
\begin{equation*}
    \langle 2,0,a \rangle = -\langle 1,1,a \rangle +\frac{1}{2}\langle 0,0,a+2\rangle.
\end{equation*}

\textbf{Inductive step:} Assume that for any $k\in \{ 1,2,\dots, n-1 \}$,we have 
\begin{align*}
    \langle n-1,n-1-k,a \rangle = 
        \sum_{j=0}^{\lfloor k/2 \rfloor} c_j(-1)^{k+j} \langle n-1-k+j, n-1-k+j, a+k-2j \rangle
\end{align*}
and let us prove for $k\in \{ 1,2\dots, n\}$ that 
\begin{align*}
    \langle n,n-k,a \rangle = 
        \sum_{j=0}^{\lfloor k/2 \rfloor} c_j(-1)^{k+j} \langle n-k+j, n-k+j, a+k-2j \rangle.
\end{align*}

\noindent Namely, if $k=1$ we have
\begin{align*}
    \langle n,n-1,a\rangle = -\frac{1}{2}\langle n-1, n-1, a+1 \rangle.
\end{align*}
If $k=2$, we have
\begin{align*}
    \langle n,n-2,a\rangle &= -\langle n-1,n-1,a\rangle-\langle n-1,n-2,a+1\rangle\\
    &=-\langle n-1,n-1,a\rangle+\frac{1}{2}\langle n-2,n-2,a+2\rangle.
\end{align*}
Now in case $k\in \{ 3,4,\dots ,n \}:$
\begin{equation*}
    \langle n,n-k,a\rangle=-\langle n-1,n-k+1,a\rangle -\langle n-1,n-k,a+1\rangle.
\end{equation*}

By inductive hypothesis, 
\begin{align}
    -\langle n-1,n-k+1,a\rangle &
    = -\sum_{j=0}^{\lfloor \frac{k-2}{2}\rfloor}c_j (-1)^{k-2+j}\langle n-1-k+2+j,n-1-k+2+j, a+k-2-2j \rangle \nonumber\\
    &= -\sum_{j=0}^{\lfloor \frac{k-2}{2}\rfloor}c_j (-1)^{k-2+j} \langle n+1+j-k, n+1+j-k, a+k-2(j+1) \rangle \nonumber\\
    &=-\sum_{j=1}^{\lfloor k/2 \rfloor} c_{j-1}(-1)^{k-3+j}\langle n+j-k,n+j-k, a+k-2j \rangle. \label{6}
\end{align}

Also, 
\begin{equation}
    -\langle n-1, n-k, a+1 \rangle =  -\sum_{j=0}^{\lfloor \frac{k-1}{2}\rfloor}c_j (-1)^{k-1+j} \langle n+j-k,n+j-k, a+k-2j \rangle .\label{7}
\end{equation}
Thus, from \eqref{6} and \eqref{7} we conclude 
\begin{align*}
    \langle n, n-k, a \rangle 
    &= \sum_{j=1}^{\lfloor k/2 \rfloor} c_{j-1} (-1)^{k-2+j}\langle n+j-k, n+j-k, a+k-2j \rangle \\& \hspace{5mm} +\sum_{j=0}^{\lfloor \frac{k-1}{2}\rfloor}c_j (-1)^{k+j} \langle n+j-k, n+j-k, a+k-2j \rangle \\
    & = c_0(-1)^k\langle n-k, n-k, a+k \rangle \\
    & \hspace{5mm} + \sum_{j=1}^{\lfloor k/2 \rfloor} c_j^* (-1)^{k+j}\langle n+j-k, n+j-k, a+k-2j \rangle\\
    & = \sum_{j=0}^{\lfloor k/2 \rfloor} c_j^* (-1)^{k+j}\langle n+j-k, n+j-k, a+k-2j \rangle,
\end{align*}
where $c_j^* = c_{j-1}+c_j$ except maybe for $j=0$ or $j=\lfloor k/2\rfloor$.

\end{proof}

\begin{corollary}
    \label{evenodd}
    For for $u$ and $\rho$ sufficiently well-behaved and any $m\in \Z^+$ we have
    \begin{align}
        \langle 2m,0,0 \rangle &= \sum_{j=0}^m c_j(-1)^j \langle j,j,2(m-j)\rangle \label{8}\\
        \langle 2m+1,0,0 \rangle &= \sum_{j=0}^m c_j(-1)^{j+1} \langle j,j,1+2(m-j)\rangle. \label{9}
    \end{align}
\end{corollary}

\section{Local theory results}\label{sec3}

We then introduce an alternative local well-posedness result in $H^s$, building upon the refined smoothing effect presented in the previous section.

\begin{theorem}[Local theory]
    \label{LWPHS}
Consider the IVP \eqref{P} where $p>5/2$ and $u_0\in H^s(\mathbb{R})$ with $s\geq 0$. Then there exist $T_0>0$ and a unique solution $u\in C([0,T_0]; H^s(\mathbb{R}))$ in the class defined by
\begin{equation}
    \label{11}
    \| D_x^s\partial_xu \|_{L_{T_0}^2L_x^4} + \| \partial_xu \|_{L_{T_0}^2L_x^4} + \| u\|_{L_{T_0}^2L_x^4}+\| D^su \|_{L_{T_0}^2L_x^4} <\infty .
\end{equation}
Moreover, for any $T'\in (0,T_0)$ there is a neighborhood $V$ of $u_0$ in $H^s(\R)$ such that the data-solution map $u_0\mapsto u$ from $V$ into the class defined by \eqref{11} with $T'$ instead of $T_0$ is Lipschitz.
\end{theorem}

\begin{proof}
    Let us first define some quantities. For $T>0$, consider
    \begin{equation*}
        A(a,b,s):=A(a,b,s)(T):= e^{\eta T}T^{1/a}+(2\alpha (a,b,s))^{-1/a}T^{\alpha(a,b,s)}.
    \end{equation*}

Define $A_2:= A(2,4,0)$, $A_3:=A(2,4,s)$ and $A_4:=A(2,4,1)$.

For simplicity, we present the details when $s \leq 1$. The general case follows by considering that derivatives fit linearly into the estimates.

For $T,z>0$ to be defined, set the following norms:
\begin{equation*}
    \begin{split}
        &\lambda_1^T(w):=\| w\|_{L_T^\infty H^s},\hspace{3mm}\lambda_2^T(w):= A_2^{-1}\| w\|_{L_T^2L_x^4}, \hspace{3mm} \lambda_3^T(w):=A_3^{-1}\| D_x^s w\|_{L_T^2L_x^4},\\ &\lambda_4^T(w):= \| D_x^s\partial_x w\|_{L_T^2L_x^4} \mbox{   and   } \lambda_5^T(w):= \| \partial_xw\|_{L_T^2L_x^4};
    \end{split}
\end{equation*}
and consider the ball
\begin{equation*}
    X_z^T:=\{ w\in C([0,T];H^s(\mathbb{R})) \mid \Lambda^T(w):=\sum_{i=1}^5\lambda_i^T(w)\leq z \}.
\end{equation*}

Let us prove that $\Psi_{u_0}:X_z^T\to X_z^T$ defined by
\begin{equation}
    \label{inteq}
    \Psi_{u_0}(w)(x,t)=V(t)u_0(x)-\int_0^t V(t-t')(u\partial_xu)(x,t')dt'
\end{equation}
 defines a contraction.

Let $u\in X_a^T$. For $t\in [0,T]$:
\begin{equation}
    \label{H1}
    \begin{split}
        \| \Psi(u)(t) \|_{L_x^2} 
    & \leq \| V(t)u_0 \|_{L_x^2}+ \int_0^T \| V(t-t')(u\partial_x u)(t') \|_{L_x^2}dt' \\
    & \leq e^{\eta T} \left\{ \| u_0\|_{L_x^2} +\int_0^T \| (u\partial_x u)(t') \|_{L_x^2}dt' \right\}\\
    & \leq e^{\eta T} \{ \| u_0\|_{L_x^2} + \| u \|_{L_T^2L_x^4}\| \partial_x u \|_{L_T^2L_x^4} \}\\
    & \leq e^{\eta T} \{ \| u_0 \|_{L_x^2}+A_2\Lambda^T(u)^2 \}.
    \end{split}
\end{equation}

Now, \begin{equation}
    \label{H2}
    \begin{split}
        \|D^s_x\Psi(u)(t)\|_{L^2_x}&\le \|D_x^sV(t)u_0\|_{L^2_x}+\int_0^T\|D^s_xV(t-t')(u\partial_x u)(t')\|_{L^2_x}dt'\\&\le e^{\eta T}\left\{\|D_x^su_0\|_{L^2_x}+\int_0^T\|D_x^s(u\partial_x u)(t')\|_{L^2_x}dt' \right\}.
    \end{split}
\end{equation}

Using Theorem \ref{graf}, for $t'\in [0,T]$ we have
\begin{equation}
    \label{H3b}
    \| D_x^s(u\partial_xu)(t') \|_{L_x^2}\leq c \| u(t') \|_{L_x^4}\| D_x^s\partial_xu(t') \|_{L_x^4}+c\| D_x^s u(t') \|_{L_x^4}\| \partial_xu(t') \|_{L_x^4}.
\end{equation}
Then,
\begin{equation}
    \label{dosestre}
    \begin{split}
        \int_0^T \| D_x^s(u\partial_x u)(t') \|_{L_x^2}dt' 
    &\leq c \| u \|_{L_T^2L_x^4}\| D_x^s \partial_xu \|_{L_T^2L_x^4} +c\| D_x^su \|_{L_T^2L_x^4} \| \partial_xu \|_{L_T^2L_x^4}\\
    &\leq cA_2\lambda_2^T(u)\lambda_4^T(u)+cA_3\lambda_3^T(u)\lambda_5^T(u).
    \end{split}
\end{equation}

Hence,
\begin{equation}   
    \| D_x^s\Psi(u)(t) \|_{L_x^2}\leq ce^{\eta T} \{ 
 \| u_0 \|_{H^s}+(A_2+A_3)\Lambda^T(u)^2  \}. \label{H3}
\end{equation}

From \eqref{H1} and \eqref{H3} we obtain
\begin{equation}
    \lambda_1^T(\Psi(u))\leq ce^{\eta T}\| u_0 \|_{H^s} +ce^{\eta T}(A_2+A_3)\Lambda^T(u)\label{H4}
\end{equation}

For $\lambda_2^T$, we use corollary \ref{cor2} with $s=0$ as follows:
\begin{equation}
    \label{H5}
    \begin{split}
         \lambda_2^T(\Psi(u))&\leq A_2^{-1} \left\{ \| V(t)u_0 \|_{L_T^2L_x^4} +\int_0^T \|V(t-t')(u\partial_xu)(t') \|_{L_T^2L_x^4}dt'  \right\}\\
    &\leq c(\eta,M)\left\{ \|u_0 \|_{L_x^2} +\int_0^T \| u\partial_x u(t')\|_{L_x^2}dt' \right\}\\
    & \leq c(\eta,M)\left\{ \| u_0 \|_{L_x^2}+ A_2\lambda_2^T(u)\lambda_5^T(u) \right\}\\
    &\leq c(\eta, M)\{ \| u_0\|_{L_x^2}+A_2\Lambda^T(u)^2 \}
    \end{split}
\end{equation}

For $\lambda_3^T$ we use corollary \ref{cor2} with $b=4$ to get 
\begin{equation}
    \label{H6}
    \begin{split}
        \lambda_3^T(\Psi(u))
    &\leq A_3^{-1} \left\{ \| D_x^s V(t)u_0 \|_{L_T^2L_x^4} +\int_0^T \| D_x^sV(t-t')(u\partial_x u)(t') \|_{L_T^2L_x^4}dt'  \right\}\\
    & \leq c(\eta, M)\left\{ \| u_0 \|_{L^2} +\int_0^T \| u\partial_x u(t') \|_{L_x^2}dt' \right\} \\
    &\leq c(\eta, M) \{ \| u_0 \|_{L_x^2} +A_2(\Lambda^T(u))^2 \}.
    \end{split}
\end{equation}

In what comes to $\lambda_4^T$, we use corollary \ref{cor2} with $s=1$ to get
\begin{equation}
    \label{H7}
    \begin{split}
        \lambda_4^T(\Psi(u)) 
    &\leq \| D_x^s \partial_xV(t)u \|_{L_T^2L_x^4} + \int_0^T \| D_x^s \partial_x V(t-t')(u\partial_x u)(t') \|_{L_T^2L_x^4}dt'\\
    &\leq \| D_x^1 V(t)\mathcal{H}D_x^su \|_{L_T^2L_x^4}+\int_0^T \| D_x^1 V(t-t')\mathcal{H}D_x^s(u\partial_x u)(t') \|_{L_T^2L_x^2}dt'\\
    &\leq c(\eta,M)A_4 \{ \| \mathcal{H}D_x^s u \|_{L_x^2} +\int_0^T \| \mathcal{H}D_x^s(u\partial_x u)(t') \|_{L_x^2}dt'  \} \\
    &\leq c(\eta,M)A_4 \{ \| u_0 \|_{H^s} +(A_2+A_3)\Lambda^T(u)^2 \},
    \end{split}
\end{equation}

where in the last inequality we used \eqref{dosestre}.

Finally, for $\lambda_5^T(\Psi(u))$ we combine the ideas used to get \eqref{H1} and \eqref{H7} to obtain
\begin{equation}
    \label{H8}
    \begin{split}
        \lambda_5^T(\Psi(u)) & \leq \| D_x^1 V(t)\mathcal{H}u_0 \|_{L_T^2L_x^4} +\int_0^T \| D_x^1V(t-t')\mathcal{H}(u\partial_x u)(t') \|_{L_T^2L_x^4}dt' \\
    &\leq c(\eta,M)A_4 \{ \| u_0 \|_{L_x^2} +A_2\Lambda^T(u)^2 \}.
    \end{split}
\end{equation}

Gathering estimates \eqref{H1}-\eqref{H8} we obtain 
\begin{align}
    \Lambda^T(\Psi(u)) &\leq c(\eta,M) \{ (e^{\eta T}+A_4)\| u_0 \|_{H^s} +(e^{\eta T}+1+A_4)(A_2+A_3)\Lambda^T(u)^2  \}\\
    &\leq c(\eta,M)(e^{\eta T}+A_4+1)\{ \| 
u_0 \|_{H^s}+(A_2+A_3)\Lambda^T(u)^2 \}.\label{hs}
\end{align}

Note that if $T\leq 1$ then $e^{\eta T}+A_4+1\leq c(e^\eta +1)$. Set $z_0:=2c^*\| u_0 \|_{H^s}$ where $c^*=c(\eta,M )c(e^{\eta}+1).$ There exists $\tilde{T}>0$ such that $(A_2+A_3)(\tilde{T})<(2c^*z_0)^{-1}$. Let $T_0:=\min(1,\tilde{T})$. We conclude $\Lambda^T(\Psi(u))\leq z_0$. 

Following standard-by-now arguments it can be seen $\Psi$ defines a contraction in $X_{z_0}^{T_0}$, completing the proof.
\end{proof}

\begin{remark}
    \label{another}
    Note that if $u\in C([0,T_0];H^{0^+}(\R))$ is the solution provided by Theorem \ref{LWPHS}, using Proposition \ref{infinity} we have for $T\le T_0$ that
    \begin{equation}
        \begin{split}
            \|D^{\frac{p}{2}}_xu\|_{L^\infty_xL^2_{T}}&\le \|D^{{\frac{p}{2}}^-}_xV(t)D^{0^+}u_0(x)\|_{L^\infty_xL^2_T}+\int_0^T\|D^{{\frac{p}{2}}^-}_xV(t-t')D^{0^+}(u\partial_x u)(t')\|_{L^\infty_xL^2_T}dt'\\&\le c(\eta,M,p,T)\left\{\|D^{0^+}u_0\|_{L^2_x}+\int_{0}^T\|D^{0^+}(u\partial_x u)(t')\|_{L^2_x}dt' \right\}\\&\le c(\eta,M,p,T_0,\Lambda^{T_0}(u)).
        \end{split}
    \end{equation}
\end{remark}

\begin{theorem}[Weighted local theory]
    \label{wLWP} Consider the IVP \eqref{P} where $p>3$ and $\Phi_1(\xi)$ satisfies the conditions \eqref{AA} and \eqref{BB} with $g(x)\leq c(1+ |x|^m)$, $m<p$. Assume $u_0\in H^s(\mathbb{R})\cap L^2(|x|^{2r}dx)$ with $s>0$ and $0\leq r\leq \frac{s}{p-1}$. Then, there is $T>0$ and a unique solution $u\in C([0,T],Z_{s,r})$ in the class defined in \eqref{11} with $T$ instead of $T_0$ satisfying also \begin{equation}
        \label{also}\| \partial_x u \|_{L_T^2L_x^\infty}<\infty.
    \end{equation} Moreover, for any $T'\in (0,T)$ there is a neighborhood $V$ of $u_0$ in $Z_{s,r}$ such that the data-solution map $u_0\mapsto u$ from $V$ into the class defined by \eqref{11} and \eqref{also} with $T'$ instead of $T_0$ or $T$ is Lipschitz.
\end{theorem}

\begin{proof}
    For $T>0$, consider $A_6:=A(2,\infty ,1)(T)$. Define the norms $\lambda_6^T(w):= A_6^{-1}\| \partial_x u \|_{L_T^2L_x^\infty}$ and $\lambda_7^T(w):=\| |x|^ru \|_{L_T^\infty L_x^2}$. For $z>0$ to be defined, set $Y_z^T:=\{ u\in C([0,T];Z_{s,r}) \mid \Omega^T(u)\le z \}$ where $\Omega^T(u):=\Lambda^T(u)+\lambda_6^T(u)+\lambda_7^T(u)$. Let us see that the integral operator $\Psi_{u_0}:Y_z^T\to Y_z^T$ defines a contraction for some $T,z>0$. Let $u\in Y_z^T$. Using Corollary \ref{cor1} with $s=1$, $a=2$ we get 
    \begin{equation}
        \label{W1}
        \begin{split}
            \lambda_6^T(\Psi(u)) & = A_6^{-1} \| D_x^1 \mathcal{H}\Psi(u) \|_{L_T^2L_x^\infty} \\
        & \leq A_6^{-1} \| D_x^1 V(t)\mathcal{H}u_0 \|_{L_T^2L_x^\infty} +\int_0^T \| D_x^1 V(t-t')\mathcal{H}(u\partial_x u)(t') \|_{L_T^2L_x^\infty}dt'\\
        &\leq c(\eta, M) \{  \| \mathcal{H}u_0 \|_{L_x^2} +\int_0^T \| \mathcal{H}(u\partial_x u)(t') \|_{L_x^2}dt' \}\\
        &\leq c(\eta,M) \{ \| u_0 \|_{L_x^2}+A_2\Lambda^T(u)^2   \}.
        \end{split}
    \end{equation}

    On the other hand, according to Theorem \ref{troca}, for $t\in [0,T]$ 
    \begin{equation}
        \| |x|^rV(t)u_0 \|_{L_x^2}\leq c(1+T) \| u_0 \|_{H^s}+c\| |x|^ru_0 \|_{L_x^2}.\label{W1n}
    \end{equation}

Similarly,
\begin{align}
    &\int_0^T \| |x|^r V(t-t')(u\partial_x u)(t')  \|_{L_x^2}dt' \nonumber\\
    &\hspace{5mm}\leq c(1+T) \int_0^T \{ \| (u\partial_x u)(t') \|_{H^s} +\| |x|^r(u\partial_x u)(t') \|_{L_x^2} \}dt'\nonumber  \\    
    & \hspace{5mm}\leq c(1+T)\int_0^T \| (u\partial_x u )(t') \|_{H^s}dt' +cT^{1/2}\| |x|^ru \|_{L_T^\infty L_x^2}  \| \partial_x u \|_{L_T^2 L_x^\infty} \nonumber\\
    &\hspace{5mm}\leq c(1+T)(A_2+A_3)\Lambda^T(u)^2+cT^{1/2}A_6\lambda_7^T(u)\lambda_6^T(u) \label{W2b} \\
    &\hspace{5mm}\leq c(1+T)(A_2+A_3+T^{1/2}A_6)\Omega^T(u)^2.\label{W2} 
\end{align}

    Thus, combining \eqref{W1n} and \eqref{W2} we conclude 
    \begin{equation}
        \label{W3}
        \lambda_7^T(\Psi(u)) \leq c(1+T) \{  \| u_0 \|_{Z{s,r}}+(A_2+A_3+T^{1/2}A_6) \Omega^T(u)^T \}.
    \end{equation}

From \eqref{hs}, \eqref{W1} and \eqref{W3}  it follows that 
\begin{equation*}
    \Omega^T(\Psi(u))\leq c(\eta ,M) (e^{\eta T}+A_4+T+2) \{ \| u_0 \|_{Z_{s,r}} +(A_2+A_3+T^{1/2}A_6)\Omega^T(u)^2 \}.
\end{equation*}
Note that if $T\le 1$ then $e^{\eta T}+A_4+T+2\le c(e^n+1)+2$. Set $z:=2c^{**}\| u_0 \|_{Z_{s,r}}$ where $c^{**}:=c(\eta,M)(c(e^\eta+1)+2)$. There exists $\tilde{T}>0$ such that $(A_2+A_3+\tilde{T}^{1/2}A_6)(\tilde{T})<(2c^{**}z)^{-1}$. Set $T:=\min \{ 1,\tilde{T}\}$. We include $\Omega^T(\Psi(u))\leq z$.

The rest of the proof follows from standard-by-now arguments.
\end{proof}

Theorem \ref{wLWP} provides an existence time $T$ that is shorter than the existence time in $H^s(\R)$, as given by Theorem \ref{LWPHS}. Nevertheless, iterations of Theorem \ref{wLWP} allow the solution to extend to the same time $T_0$ as in $H^s(\R)$ and can be estate as a persistence result.
We now enunciate the persistence result associated to Theorem \ref{wLWP}.
\begin{theorem}[Persistence]\label{sametime}
    Assume $s>0$, $p>3$ and let $u\in C\left([0,T_0];H^s(\R) \right)$ be the solution provided by Theorem \ref{LWPHS}. Suppose there exists $t^*\in[0,T_0]$ such that $|x|^r u(t^*)\in L^2(\R)$ with $0<r\le \frac{s}{p-1}$. Then, $u\in C\left([0,T_0];Z_{s,r} \right)$.
\end{theorem}
\begin{proof}
    Recall from the proof of Theorem \ref{LWPHS} that $T_0\le 1$ is taken satisfying 
\begin{equation*}
    (A_2+A_3)(T_0)<\frac{1}{2c^*z_0}, \mbox{ where } c^*=c(\eta,M)c(e^\eta+1) \mbox{ and } z_0=2c^*\|u(0)\|_{H^s}.
\end{equation*}

We first establish some \textit{a priori} bounds for $T\le T_0$. From \eqref{W1} we have, 
\begin{equation}\label{ST1}
    \begin{split}
        \lambda_6^T(u)&\le c(\eta,M)\left\{\|u(0)\|_{L^2_x}+A_2(T)\Lambda^T(u)^2\right\}\\
        &\le c(\eta,M)\left\{\|u(0)\|_{L^2_x}+A_2(T_0) z_0^2\right\}\\
        &\le c(\eta,M)\left\{\|u(0)\|_{L^2_x}+(2c^*z_0)^{-1}z_0^2\right\}\\&
        \le c(\eta,M)\left\{\|u(0)\|_{H^s}\right\}
    \end{split}
\end{equation}
Similarly, according to \eqref{W2b} and using \eqref{ST1}, for $T\le T_0$ we have
\begin{align}
        \lambda^T_7(u)&\le c(1+T)\|u(0)\|_{H^s}+c\||x|^r u\|_{L^2_x}+(A_2+A_3)(T)\Lambda^T(u)^2+cT^{1/2}A_6(T)\lambda_7^T(u)\lambda_6^{T}(u)\nonumber \\
        &\le 2c\left\{\|u(0)\|_{Z_{s,r}}+(A_2+A_3)(T_0)z_0^2+T_0^{1/2}A_6(T_0)\lambda_7^T(u)\lambda_6^{T}(u) \right\}\nonumber\\
        &\le 2c\left\{\|u(0)\|_{Z_{s,r}}+(2c^*z_0)^{-1}z_0^2+c(\eta,M)T_0^{1/2}A_6(T_0)\|u(0)\|_{H^s}\lambda_7^T(u)  \right\} 
        \nonumber\\\label{ST2}
        &\le 2c\left\{\|u(0)\|_{Z_{s,r}}+c(\eta,M)T_0^{1/2}A_6(T_0)\|u(0)\|_{H^s}\lambda_7^T(u)  \right\}.
\end{align}

Since $T_0\le 1$ and $\alpha(2,\infty,1)\le\alpha(2,4,0)$ we have $A_6(T_0)\le A_2(T_0)$. Therefore
\begin{equation*}
    \begin{split}
        c(\eta,M)T_0^{1/2}A_6(T_0)\|u(0)\|_{H^s} &\le c(\eta,M)A_2(T_0)\|u(0)\|_{H^s} \\ &\le c(\eta,M)(2c^*z_0)^{-1}\|u(0)\|_{H^s}\\&\le c(\eta,M)\left( 4c(\eta,M)^2\|u(0)\|_{H^s}c^2(e^\eta+1)^2\right)^{-1}\|u(0)\|_{H^s}\\&\le \left(4c(\eta,M)c^2(e^\eta+1)^2 \right)^{-1}=:B<1.
    \end{split}
\end{equation*}
Thus, from \eqref{ST2} we get for any $T\le T_0$ that
\begin{equation}
    \label{ST3}
    \lambda_7^T(u)\le (1-B)2c\|u(0)\|_{Z_{s,r}}=:c^{\#}\|u(0)\|_{Z_{s,r}}.
\end{equation}

Now, consider $\Tilde{u}_0(x):=u(t^*)$ and the IVP \eqref{P} with $\Tilde{u}_0(x)$ as initial datum. Using Theorem \ref{wLWP} we obtain a solution $\Tilde{u}\in C\left([t^*,t^*+T];Z_{s,r}\right)$. By uniqueness we have 
\begin{equation}
    u\in C\left([0,T_0];H^s(\R)\right)\cap C\left([t^*,t^*+T];Z_{s,r}\right), 
\end{equation}
where $T$ was chosen so that $(A_2+A_3+T^{1/2}A_6)(T)<(2c^{**}\Tilde{z})^{-1}$ with $\Tilde{z}=2c^{**}\|\Tilde{u}(0)\|_{Z_{s,r}}$. In particular, $T\le T_0$. 

Note \begin{equation}
    \Tilde{z}\le 2c^{**}\|\Tilde{u}(0)\|_{H^s}+2c^{**}\||x|^r\Tilde{u}(0)\|_{L^2_x}. \label{ST4}
\end{equation}

Since \begin{equation*}
    \|\Tilde{u}(0)\|_{H^s}=\|u(t^*)\|_{H^s}\le \lambda_1^T(u)\le z_0\le2c^{*}\|u(0)\|_{Z_{s,r}}
\end{equation*}
and
\begin{equation*}
    \||x|^r \Tilde{u}(0)\|_{L^2_x}=\||x|^r u(t^*)\|_{L^2_x}\le \lambda_7^T(u)\le c^{\#}\|u(0)\|_{Z_{s,r}},
\end{equation*}
 we conclude from \eqref{ST4} that 
 \begin{equation*}
     \Tilde{z}\le 4c^{*}c^{**}\|u(0)\|_{Z_{s,r}}+2c^{**}c^{\#}\|u(0)\|_{Z_{s,r}}=:z^\#.
 \end{equation*}

Fix $T^\#>0$ such that $(A_2+A_3+{(T^\#)}^{1/2}A_6)(T^\#)<(4cz^\#)^{-1}$ (in case $T^\#>1$, set $T^\#$ equal to 1). We have $0<T^\#\le T$. Moreover, since estimates \eqref{ST1} and \eqref{ST2} are valid over the whole interval $[0,T_0]$, for any $T$ coming from iterations of Theorem \ref{wLWP} we also have $T^\#\le T$; which allow us to reapply the weighted local well-posedness theory $T_0/T^\#$-times to obtain $u\in C\left([0,T_0];Z_{s,r} \right)$. \end{proof}

\section{Regularity vs. decay}\label{sec4}
In this section, we investigate the optimality of the relationship between the regularity and decay indices, denoted as $s$ and $r$ respectively. 

\begin{theorem}[Optmality]
    \label{optimal}
    Consider the IVP \eqref{P} with $\Phi(\xi)=-|\xi|^{2k}+\xi^{2k-3}|\xi|^2$ where $k\in\{2,3,\dots\}$; \textit{i.e.}, with $Lu=(-1)^k\partial_x^{2k}u-i\partial_x^{2k-1}u$. Let $u_0\in L^2(\R)$ and assume the solution $u\in C\left([0,T_0];L^2(\R)\right)$ provided by Theorem \ref{LWPHS} is so that for two different times $t_0, t_1$ and $\alpha>0$
    \begin{equation}\label{part}
        |x|^\alpha u(t_i) \in L^2(\R), \ \ i=0,1.
    \end{equation} 
    Also assume that if $\alpha\in(\frac{n-1}{2},\frac{n}{2}]$ for $n\in \mathbb{Z}^+$, there is a time $t_*\in[0,T_0]$ such that \begin{equation}
        \label{extra}\|D^{((n-1)k)^+}_x u(t_*)\|_{L^2_x}<\infty.
    \end{equation}
    Then, $u\in C\left([0,T_0];H^{2k\alpha}(\R)\right)$.
\end{theorem}
\begin{remark}
    The slight additional regularity required in hypothesis \eqref{extra} appears to be essential for completing the approximation argument developed in the proof of Theorem \ref{optimal}. It is important to note that this requirement always involves only a small increment in regularity. Each step in the proof of Theorem \ref{optimal} necessitates an additional $\varepsilon$ of regularity. For instance, when $\alpha\in(0,1/2]$, from $H^{0^+}$ the regularity is improving up to $H^{2k\alpha}$ but if $\alpha\in(1/2,1]$, in particular, hypothesis \eqref{part} holds for $\alpha=1/2$ and therefore we can conclude (from $H^{0^+}$) that $u(t)$ is already in $H^k$ which transform hypothesis \eqref{extra} for $\alpha\in(1/2,1]$ into a requirement of extra $\varepsilon$ of regularity. 
\end{remark}
The proof of Theorem \ref{optimal} follows using induction and approximation arguments as done in \cite[Theorem 1.4]{ILP2013} and \cite[Theorem 1.1]{MP2}. We present here a simplified proof omitting some technicalities. 

\begin{proof}[Proof of Theorem \ref{optimal}]
We argue by cases depending on the size of $\alpha$. 

\textbf{Case $\alpha\in(0,1/2]$.}\\  
Let $\rho_N=\rho$ be a regularization of \begin{equation}
    \rho_N(x)=\begin{cases}
        \langle x\rangle^{2\alpha}, &x\in[-N,N]\\ 0 & x\in[-2N,2N]^c
    \end{cases}
\end{equation}
 satisfying $|\rho^{(j)}|\le c$, $j=1,2,\dots$, with $c$ independent on $N$, and $0<\rho_N(x)\le c\langle x\rangle^{2\alpha}$ for all $N\in \mathbb{N}$. 

 Without loss of generality assume $t_0=0$. 

 Multiplying the equation by $(-1)^k u \rho$ and integrating by parts in the spatial variable using formulae \eqref{8} and \eqref{9} we get

\begin{equation}\label{o1}
    \begin{split}
        (-1)^k&\left\{\frac{1}{2}\frac{d}{dt}\int u^2\rho dx +\frac{3}{2}\int (\partial_x u)^2\rho'dx -\frac{1}{2}\int u^2 \rho^{(3)}dx -\frac{1}{3}\int u^3 \rho'dx\right.\hspace{35mm} \\&\hspace{55mm}+\eta (-1)^k\sum_{j=0}^{k}c_j (-1)^j \int (\partial_x^ju)^2 \rho^{(2(k-j))}dx\\&\hspace{55mm}\left.-\eta i \sum_{j=0}^{k-1}c_j (-1)^{j+1} \int (\partial_x^ju)^2 \rho^{1+(2(k-1-j))}dx\right\}=0.
    \end{split}
\end{equation}
 We can rewrite \eqref{o1} as follows
 \begin{equation}\label{o2}
     \begin{split}
         \eta c_k\int (\partial_x^k u)^2\rho
dx=(-1)^k&\left\{-\frac{1}{2}\frac{d}{dt}\int u^2\rho dx -\frac{3}{2}\int (\partial_x u)^2\rho'dx +\frac{1}{2}\int u^2 \rho^{(3)}dx +\frac{1}{3}\int u^3 \rho'dx \right. \\ &\left.-\eta(-1)^k\sum_{j=0}^{k-1}c_j(-1)^j\int (\partial_x^j u)^2\left(\rho^{(2(k-j))}+i\rho^{(1+2(k-1-j))}\right)dx\right\}.  \end{split}
 \end{equation}

For $0\le j\le k-1<k$, using Corollary \ref{cor4} with $s=j$ we have
\begin{equation}\label{o2n}
    \begin{split}
        \left(\int_0^{t_1} \int (\partial_x^j u)^2dxdt \right)^{1/2}&\le \|\partial_x^ju\|_{L^2_{xT_0}}\\&\le \|\partial_x^j V(t)u_0\|_{L^2_{xT_0}}+\int_{0}^{T_0} \|\partial_x^j V(t-t')u\partial_x u (t')\|_{L^2_{xT_0}}dt'\\&\le c(\eta,M,T_0)\left\{\|u_0\|_{L^2_x}+\int_0^{T_0}\|u\partial_x u(t')\|_{L^2_x}dt' \right\}\\& \le c(\eta,M,T_0, \Lambda^{T_0}(u)). 
     \end{split}
\end{equation}

Similarly, using Corollary \ref{cor1} with $s=0$, it can be seen that
\begin{equation}
    \|u\|_{L^2_TL^\infty_x}\le c(\eta, M, T_0, \Lambda^{T_0}(u)). \label{o3b}
\end{equation}

Since $|\rho^{(j)}|\le c$ for $j=1,2,\dots$, integrating equation \eqref{o2} over $[0,t_1]$, using \eqref{o2n}, \eqref{o3b} and the hypothesis on $u_0$ and $u(t_1)$ we get
\begin{equation}\label{o3}
    \int_0^{t_1}\int (\partial_x^ku)^2\rho_Ndxdt\le c,
\end{equation}
where $c$ depend on $\|\langle x\rangle^\alpha u_0 \|_{L^2_x}$, $\|\langle x\rangle^\alpha u(t_1) \|_{L^2_x}$, $T_0$ and $\Lambda^{T_0}(u)$. 
(It is worth mentioning that when dealing with the approximation argument, \eqref{o3} is obtained by taking into account Remark \ref{another} and a similar procedure as done in (3.8) in \cite{MP2} together with the continuous dependence on $\Lambda^{T_0}$). 

Since $c$ in \eqref{o3} does not depend on $N$, we have
\begin{equation}\label{o4}
    \int_0^{t_1}\int (\partial_x^k u)^2 \langle x\rangle^{2\alpha}dxdt \le \liminf_{N\to \infty}\int_0^{t_1}\int (\partial_x^k u)^2\rho_N dxdt<\infty.
\end{equation}

In particular $\langle x\rangle^{\alpha-1/2}\partial_x^k u (t)\in L^2(\R),$ \textit{a.e. }$t\in[0,t_1]$. Moreover, using that $\alpha-1/2\le 0$, \eqref{o2n} and the Leibniz rule it can be seen that 
\begin{equation}
    \label{o5}J^k(\langle x\rangle^{\alpha-1/2}u)\in L^2(\R)\ \mbox{\textit{a.e. }}t\in[0,t_1].
\end{equation}

Taking into account \eqref{o4} and isolating $\frac{d}{dt}\int u^2\rho dx$ instead of $\int (\partial_x^k u)^2\rho dx$ in \eqref{o1}, it can be seen, in a similar fashion as done to obtain \eqref{o4}, that for any $t\in(0,t_1)$ we have
\begin{equation}
    \langle x\rangle^{\alpha}u(t)\in L^2(\R). \label{o6}
\end{equation}
Let $t^*\in[0,t_1]$ be such that $\langle x\rangle^{\alpha}u(t^*)$ and $J^k(\langle x\rangle^{\alpha-1/2}u(t^*))$ are in $L^2(\R)$. set $f=\langle x\rangle^{\alpha-1/2}u(t^*)$. From \eqref{o5} and \eqref{o6}, interpolating using Lemma \ref{interp} with $\theta=2\alpha$, $a=k$ and $b=1/2$, we get
\begin{equation}
    \label{o6n}
    \|J^{2k\alpha}u(t^*)\|_{L^2_x}=\|J^{\theta k}\left(\langle x\rangle^{(1-\theta)/2}f \right)\|_{L^2_x}\le c\|J^kf\|_{L^2_x}^\theta \|\langle x\rangle^{1/2}f\|_{L^2_x}^{1-\theta}<\infty.
\end{equation}
We conclude $u(t^*)\in H^{2k\alpha}(\R)$, which implies $u\in C\left([0,T_0];H^{2k\alpha}(\R) \right)$. Moreover, taking into account Theorem \ref{sametime}, we have $u\in C([0,T_0];Z_{2k\alpha, \frac{2k\alpha}{2k-1}} )$.

\textbf{Case $\alpha>1/2$}.\\
In general, when $\alpha\in (\frac{n-1}{2},\frac{n}{2}]$ the result follows by an inductive process by taking $\partial_x^{k(n-1)}$ derivatives to the equation and adjusting the exponent in $\rho_N$ to $2\alpha-n+1$ (see the proof of Theorem 1.1 in \cite{MP2}); we carry on some details when $\alpha\in(1/2,1]$. 

Assume $\alpha\in(1/2,1]$. In particular, \eqref{part} holds for $\alpha=1/2$ and therefore, from the previous case we know 
\begin{equation}
    \label{moda}u\in C([0,T_0]; Z_{k,\frac{k}{2k-1}}) \mbox{ with } \langle x\rangle^{1/2}\partial_x^k u(t)\in L^2(\R), \mbox{ \textit{a.e } in } [0,t_1]
\end{equation}

\underline{\textbf{Claim:}} If \begin{equation}
    \label{o7b} \|\langle x\rangle^{\alpha-1/2}\partial_x^{k-1}u\|_{L^2_{xt_1}}\le c(\eta,M,T_0,\Omega^{T_0}(u))
\end{equation}
holds for some $\alpha\in (1/2,1]$, then $u\in C([0,T_0];Z_{2k\alpha,\frac{2k\alpha}{2k-1}})$.

\underline{\textit{Proof of the claim:}}
Under these hypotheses we first prove \eqref{o4} remains valid when $\alpha\in (1/2,1]$. To do so, the same steps \eqref{o1}-\eqref{o4} are done, except for the fact that $|\rho'|$ is no longer bounded above by a constant independent on $N$. We provide the respective modifications on the estimates of the terms in \eqref{o2} involving $\rho'$. 

Note 
\begin{equation*}
  -\frac{3}{2}\int (\partial_x u)^2\rho'dxdt  = \frac{3}{2}\int u\partial_x^2u\rho'dx+\frac{3}{2}\int u\partial_xu\rho''dx,
\end{equation*}
with 
\begin{equation*}
    \left|\frac{3}{2}\int_0^{t_1}\int u\partial_xu\rho''dxdt \right|\le c\|u\|_{L^2_{xT_0}}\|\partial_xu\|_{L^2_{xT_0}}\le c(\eta,M,T_0,\Lambda^{T_0}(u))
\end{equation*}
and
\begin{equation*}
\begin{split}
    \left|\frac{3}{2}\int_0^{t_1}\int u\partial_x^2u\rho'dxdt \right|&\le c\|\partial_x^2 u\|_{L^2_{xT_0}}\|\langle x\rangle^{1/2}u\|_{L^2_{xt_1}}\le c\|\partial_x^2 u\|_{L^2_{xT_0}}\|\langle x\rangle^{\frac{k}{2k-1}}u\|_{L^2_{xT_0}}\\&\le c(\eta,M,T_0,\Omega^{T_0}(u)).
\end{split}
\end{equation*}

Thus, 
\begin{equation*}
    \left|-\frac{3}{2}\int_0^{t_1}\int (\partial_x u)^2\rho'dxdt\right|\le c(\eta,M,T_0,\Omega^{T_0}(u)).
\end{equation*}

Similarly, using Sobolev embeddings,
\begin{equation*}
    \left|\frac{1}{3}\int_0^{t_1}\int u^3\rho'dxdt \right|\le \|u\|_{L^\infty_{xT_0}}\|u\langle x \rangle^{1/2}u\|_{L^2_{xT_0}}\le c(\eta,M,T_0,\Omega^{T_0}(u)).
\end{equation*}

Also, according to \eqref{o7b} we have \begin{equation*}
    \left|-\eta c_{k-1}i\int_0^{t_1}\int (\partial_x^{k-1}u)^2\rho'dxdt \right|\le c\|\langle x\rangle^{\alpha-1/2}\partial_x^{k-1}u\|_{L^2_{xt_1}}\le c(\eta,M,T_0,\Omega^{T_0}(u)). 
\end{equation*}

This way, all the terms in \eqref{o2} involving $\rho'$ are bounded, which allow us to conclude 
\begin{equation}
    \label{o7bb}
    \int_0^{t_1}\int (\partial_x^ku)^2\langle x\rangle^{2\alpha}dxdt<\infty.
\end{equation}

Now, Consider $\rho_N=\rho$ as before but based on $\langle x\rangle^{2\alpha-1}$ instead of $\langle x\rangle^{2\alpha}$. Denote $w:=\partial_x^k u$. Take $\partial_x^k$ to the equation to get 
\begin{equation}\label{o7}
    \partial_t w + \partial_x^3 w + \eta L w +\partial_x^k(u\partial_xu)=0.
\end{equation}
Note the linear part of \eqref{o7} remains analogous to the original equation when seen for $w$. Multiplying \eqref{o7} by $(-1)^kw\rho$ and integrating by parts using formulae \eqref{8} and \eqref{9} we get
\begin{equation*}
     \begin{split}
         \eta c_k\int (\partial_x^k w)^2\rho
dx=(-1)^k&\left\{-\frac{1}{2}\frac{d}{dt}\int w^2\rho dx -\frac{3}{2}\int (\partial_x w)^2\rho'dx +\frac{1}{2}\int w^2 \rho^{(3)}dx\right. \\ &\left.-\eta(-1)^k\sum_{j=0}^{k-1}c_j(-1)^j\int (\partial_x^j w)^2\left(\rho^{(2(k-j))}+i\rho^{(2(k-1-j))}\right)dx\right.\\&\left.-\int \partial_x^k(u\partial_x u)w\rho dx \right\}.
  \end{split}
 \end{equation*}
Since $u\in H^k(\R)$, in terms of $w$, estimates \eqref{o2n} and \eqref{o3b} remain equal, we need to analyze the nonlinear term only. Using Leibniz rule
\begin{equation*}
    \partial_x^k (u\partial_x u)=\sum_{j=0}^{k}\binom{k}{j}\partial_x^{k-j}u\partial_x^{j}u.
\end{equation*}
This way \begin{equation*}
    -\int_0^{t_1} \int \partial_x^k(u\partial_x u)w\rho dx =-\sum_{j=0}^{k}\binom{k}{j}\int_0^{t_1}\int \partial_x^{k-j}u\partial_x^{j}u w\rho dx=:\sum_{j=0}^{k}\binom{k}{j}B_j.
\end{equation*}
We estimate the terms in three cases.

\underline{If $j=0$ or $j=k-1$:} Using Sobolev embeddings and \eqref{moda} we get
\begin{equation*}
    \begin{split}
        \left|-\int_0^{t_1}\int \partial_x^k u \partial_x u w \rho dxdt \right|&\le \int_0^{t_1}\int w^2|\partial_x u|\rho dxdt \le \|\partial_xu\|_{L^\infty_{xT_0}}\|\langle x\rangle^{\alpha-1/2}w\|_{L^2_{xt_1}}\\&\le c\|u\|_{L^\infty_{T_0}H^k}\|\langle x\rangle^{1/2}\partial_x^k u\|_{L^2_{xt_1}}\le c(\eta, M, T_0, \Lambda^{T_0}(u)).
    \end{split}
\end{equation*}

\underline{If $j=k$:} Integrating by parts,
\begin{equation*}
    -\int u\partial_x^{k+1}u w \rho dx = \int \partial_x u w^2 \rho dx = \int uw\partial_x^{k+1}u\rho dx +\int u w^2\rho'dx.
\end{equation*}
Thus, using Sobolev embedding and \eqref{moda}, we obtain
\begin{equation*}
    \begin{split}
        \left|-\int_0^{t_1}\int u \partial_x^{k+1}u w\rho dxdt \right|&\le \frac{1}{2}\int_0^{t_1}\int |\partial_xu|w^2\rho dxdt + \frac{c}{2}\int_0^{t_1}\int |u|w^2dxdt \\&\le c\|\partial_x u\|_{L^\infty_{xT_0}}\|\langle x\rangle^{1/2}w\|_{L^2_{xt_1}}+c\|u\|_{L^\infty_{xT_0}}\|w\|_{L^2_{xT_0}}\\ &\le c(\|u\|_{L^\infty_{T_0}H^k}\|\langle x\rangle^{1/2}\partial_x^k u\|_{L^2_{xt_1}}+c\|u\|^2_{L^\infty_{T_0}H^k} )\le c(\eta,M,T_0,\Lambda^{T_0}(u)). 
    \end{split}
\end{equation*}

\underline{If $j\in\{1,2,\dots,k-2\}$:} Using Sobolev embeddings and \eqref{o7bb} we get
\begin{equation*}
    \begin{split}
        \left|-\int_0^{t_1}\int \partial_x^{k-j}u\partial_x^{1+j}u w \rho dxdt \right|&\le \|\langle x\rangle^{2\alpha-1}w\|_{L^2_{xt_1}}\|\partial_x^{k-j}u\|_{L^\infty_{xT_0}}\|\partial_x^{1+j}u\|_{L^2_{xt_1}}\\&\le c(T_0) \|\langle x\rangle^{\alpha}\partial_x^k u\|_{L^2_{xt_1}}\|u\|^2_{L^\infty_{T_0}H^k}\le c(\eta, M, T_0,\Omega^{T_0}(u)). 
    \end{split}
\end{equation*}

In any case, \begin{equation}\label{falsear}
    \left|B_j\right|\le c(\eta,M,T_0,\Omega^{T_0}(u)). 
\end{equation}

Arguing as in the previous case it can be seen that for some $f=\langle x\rangle^{\alpha-1}w(t^*)$ we have $\langle x\rangle^{1/2}f$ and $J^kf$ are in $L^2(\R)$ with $t^*\in [0,t_1]$. Interpolating as in \eqref{o6n} with $\theta=2\alpha-1$, $a=k$ and $b=1/2$ we get
\begin{equation}
    \label{o8}
    \|J^{k(2\alpha-1)}w(t^*)\|_{L^2_x}< \infty.
\end{equation}
Note \eqref{o8} implies $\partial_x^k u(t^*)\in H^{2k\alpha-k}(\R)$; that is, $u(t^*)\in H^{2k\alpha}(\R)$.\\ The claim follows from the local theory in $H^s(\R)$ and Theorem \ref{sametime}. \\
\vspace{2mm}

To conclude the proof we show that \eqref{o7b} is valid for all $\alpha\in(1/2,1]$. Taking into account Lemma \ref{hilberpp} and interpolating with $\theta=\frac{k}{2k-1}$, $a=2k-1$ and $b=(\alpha-1/2)\frac{2k-1}{k}$ we have 
\begin{equation}
    \label{o9}
    \|\langle x\rangle^{\alpha-1/2}\partial_x^{k-1}u\|_{L^2_x}\le c \|D^{2k-1}_x u\|_{L^2_x}^\theta \|\langle x\rangle^b u\|^{1-\theta}_{L^2_x}.
\end{equation}

Note using Corollary \ref{cor4}, 
\begin{equation}
    \label{o10}
    \begin{split}
      \|D^{2k-1}_xu\|_{L_{xt_1}^2}&\le \|D_x^{k-1}V(t)D_x^k u_0\|_{L^2_{xt_1}}+\int_0^{T_0}\|D_x^{k-1}V(t-t')D_x^k(u\partial_x u)(t')\|_{L^2_{xt_1}}dt'\\&\le c(\eta,M,T_0)\left\{\|D_x^ku_0\|_{L^2_x} +\int_0^{T_0}\|D_x^{k}(u\partial_x u)(t')\|_{L^2_x}dt'\right\}\le c(\eta,M,T_0,\Lambda^{T_0}(u)).
    \end{split}
\end{equation}

From \eqref{o9} and \eqref{o10}, 
\begin{equation}
    \label{o11}
    \begin{split}
        \|\langle x\rangle^{\alpha-1/2}\partial_x^{k-1}u\|_{L^2_{xt_1}}&\le \left\| \|D^{2k-1}_x u\|_{L^2_x}^\theta \|\langle x\rangle^bu\|_{L^2_x}^{1-\theta} \right\|_{L^2_{t_1}}\\&\le \|D^{2k-1}_x u\|_{L^2_{xt_1}}+\|\langle x\rangle^b u\|_{L^2_{xt_1}}\\&\le c(\eta,M,T_0,\Lambda^{T_0}(u))+c(T_0)\|\langle x\rangle^b u\|_{L^\infty_{T_0}L^2_x}. 
    \end{split}
\end{equation}
It remains to estimate $\|\langle x \rangle^{b}u\|_{L^\infty_{T_0}L^2_x}$. We use a bootstrap argument as follows: note $b\le \frac{k}{2k-1}$ for $\alpha\le 3/4$, in such case accoridng to \eqref{moda}
\begin{equation}\label{o12}
    \|\langle x\rangle^b u\|_{L^\infty_{T_0}L^2_x}\le c\|\langle x\rangle^{\frac{k}{2k-1}}u\|_{L^\infty_{T_0}L^2_x}\le c(\Omega^{T_0}(u)).
\end{equation}
From \eqref{o11} and the claim, for any $\alpha\in(1/2,3/4]$ we have $u\in C([0,T_0];Z_{2k\alpha,\frac{2k\alpha}{2k-1}})$. Now, if $\alpha\in(3/4,7/8]$, we know $u\in C([0,T_0];Z_{\frac{3k}{2},\frac{3k}{4k-2}})$. Thus, under the same interpolation \eqref{o9}, instead of \eqref{o12} it is enough to have $b\le\frac{3k}{4k-2}$, which is true for $\alpha\le 7/8.$ 
In case $\alpha\in(7/8, 15/16]$, from the claim we know $u\in C([0,T_0];Z_{\frac{7k}{4},\frac{7k}{8k-4}})$ and thus, under the same interpolation \eqref{o9} it is enough to have $b\le \frac{7k}{8k-4}$, which holds for $\alpha\le\frac{15}{16}$. 

Continuing this process it can be seen that $\|\langle x\rangle^b u\|_{L^\infty_{T_0}L^2_x}$ is bounded in terms of $\Omega^{T_0}(u)$ for $\alpha\le\frac{2^n-1}{2^n}$ for any $n\in\N$. In other words, 
\begin{equation}
    \label{o13}\mbox{\eqref{o7b} is true for any }\alpha\in (1/2,1)\mbox{ with } u\in C([0,T_0],Z_{2k\alpha,\frac{2k\alpha}{2k-1}}).
\end{equation}

We finally prove \eqref{o7b} for $\alpha=1$. Since $k>1$, we have $1-\frac{3}{4k}$ belongs to $(1/2,1)$. From \eqref{o13} we know $u\in C([0,T_0];Z_{2k-3/2,\frac{2k-3/2}{2k-1}})$. Interpolating with $\theta=\frac{2k-3/2-k+1}{2k-1}=\frac{k-1/2}{2k-3/2}$, $a=2k-3/2$ and $b=\frac{1}{2\theta}=\frac{2k-3/2}{2k-1}$ we get
\begin{equation}\begin{split}
    \|\langle x\rangle^{1/2}\partial_x^{k-1}u\|_{L^2_x}&\le c\|D^{2k-3/2}_xu\|_{L^2_x}^\theta \|\langle x \rangle^{\frac{2k-3/2}{2k-1}}u\|_{L^2_x}\le c(\eta,M,T_0,\Omega^{T_0}(u)).
\end{split}
\end{equation}
\end{proof}
\begin{corollary}\label{coro}
    Let $s\in ((n-1)k,nk)$ with $n,k\in \Z^+$ and let $u\in C([0,T_0];H^s(\R))$ be the solution of \eqref{P} with $\Phi(\xi)$ as in Theorem \ref{optimal}. Suppose for two different times $t_0$, $t_1$ in $[0,T_0]$\begin{equation*}
        |x|^{\frac{s}{2k}+\alpha}u(t_i)\in L^2(\R)\ \ i=0,1; \mbox{ with }0<\alpha\le \frac{n}{2}-\frac{s}{2k}.
    \end{equation*}
    Then, $u\in C([0,T_0];H^{s+2k\alpha}(\R))$. 
\end{corollary}
\begin{proof}
    It follows directly from Theorem \ref{optimal} since $s>(n-1)k$ and $\frac{s}{2k}+\alpha \in (\frac{n-1}{2},\frac{n}{2}]$.
\end{proof}

\section{Asymmetrically weighted Sobolev spaces}\label{sec5}
For simplicity we consider the Korteweg-de Vries-Kuramoto-Sivashinsky (KdV-KS) equation 
\begin{equation}
    \label{kdvks}\begin{cases}
        \partial_t u + \partial_x^3 u +\eta(\partial_x^2+\partial_x^4)u+u\partial_xu=0, &t>0,\ x\in \R;\\ \hspace{39mm}u(0,x)=u_0(x),
    \end{cases}
\end{equation}
where $u_0\in H^s(\R)\cap L^2(e^{2bx}dx)$. Note in this case $p=4$ and $\Phi_1(\xi)=|\xi|^2$. From the proof of the results coming after, it can be seen that a similar procedure holds for $p>1$ even integer and $\Phi_1$ taken as a linear combination of terms of the form $\xi^m|\xi|^n$ with $n$ even integer and $m+n<p$. 

\begin{proposition}
    \label{trocae}
    Let $t,b>0$. If $e^{bx}f\in L^2(\R)$ we have 
    \begin{equation*}
        \|e^{bx}V(t)f\|_{L^2_x}\le c(\eta, b)e^{-t\delta}(1+e^t)\|e^{bx}f\|_{L^2_x},
    \end{equation*}
    where $V(t)$ is the group associated to the linear part of \eqref{kdvks} and $\delta=\eta(b^2+b^4)-b^3$. 
\end{proposition}
\begin{proof}
Set $w(x,t)$ to be the solution of the linear IVP \begin{equation}
    \begin{cases}
        \partial_t w+(\partial_x-b)^3w+\eta\left( (\partial_x-b)^2+(\partial_x-b)^4\right)w=0,\\ \hspace{61mm}w(0,x)=e^{bx}f(x).
    \end{cases}
\end{equation}
It can be seen that $e^{-bx}w(x,t)$ satisfies the linear IVP associated to \eqref{kdvks} and therefore $w(x,t)=e^{bx}V(t)u_0(x)$. 

Note the spatial Fourier transform of $w$ satisfies \begin{equation*}
    \widehat{w(x,t)}=\exp{\left\{-t\left\{ (i\xi-b)^3+\eta((i\xi-b)^2+(i\xi-b)^4)\right\} \right\}}\widehat{e^{bx}u_0}.
\end{equation*}
Opening the products we have
\begin{equation*}
\begin{split}
   \exp&\left\{-t\left\{(i\xi)^3-3b(i\xi)^2+3b^2(i\xi)-b^3+\eta((i\xi)^2-2b(i\xi)+b^2\right.\right.\\&\hspace{50mm}\left.\left.+(i\xi)^4-4b(i\xi)^3+6b^2(i\xi)^2-4b^3(i\xi)+b^4    )  \right\}  \right\}\\=&\exp\left\{-t(3b\xi^2-b^3+\eta(-\xi^2+b^2+\xi^4-6b^2\xi^2+b^4))\right\}\\&\hspace{50mm}\cdot\exp\left\{-it(-\xi^3+3b^2\xi+\eta(-2b\xi-4b\xi^3-4b^3\xi)) \right\}.
\end{split}
\end{equation*}
Set $\theta(\xi):=\eta\xi^4+(3b-\eta-6\eta b^2)\xi^2$ and $\mu=3b^2-2\eta b-4\eta b^3$. We have
\begin{equation}\label{hung}
\begin{split}
\|w(x,t)\|_{L^2_x}&=\|\widehat{w(x,t)}\|_{L^2_\xi}\le \|e^{-t\theta(\xi)}\|_{L^\infty_\xi}\|e^{-t\mu \partial_x}e^{-t\delta}e^{bx}f\|_{L^2_x}\\&\le\|e^{-\theta(\xi)}\|_{L^\infty_\xi}\|e^{b(x-\mu t)}e^{-t\delta}f(x-\mu t)\|_{L^2_x}
\\&\le \|e^{-t\theta(\xi)}\|_{L^\infty_\xi}e^{-t\delta}\|e^{bx}f\|_{L^2_x}.
\end{split}
\end{equation}
In case $3b-\eta -6\eta b^2>0$ we can bound directly $\|e^{-t\theta(\xi)}\|_{L^\infty_\xi}\le c$. If not, there is a constant $c_0=c(\eta, b)$ such that for all $|\xi|\ge c_0$ we have $\eta\xi^4+(3b-\eta-6\eta b^2)\xi^2\ge \frac{|\xi|^4}{2}$. Therefore, 
\begin{equation}
\begin{split}
    \|e^{-t\theta(\xi)}\|_{L^\infty_\xi}&\le \|e^{-t\theta(\xi)}\|_{L^\infty_\xi(|\xi|\le c_0)}+ \|e^{-t\theta(\xi)}\|_{L^\infty_\xi(|\xi|\ge c_0)}\\&\le \exp\left\{ -t(3b-\eta-6b^2)c_0^2\right\}+\|e^{-t\frac{|\xi|^4}{2}}\|_{L^\infty_\xi}\\&\le \exp\left\{ -t(3b-\eta-6b^2)c_0^2\right\}+1\\&\le c(\eta,b)(e^t+1).
\end{split}
\end{equation}
In any case, \begin{equation}
   \|e^{-t\theta(\xi)}\|_{L^\infty_\xi}\le c(\eta,b)(e^t+1).\label{uju}
\end{equation}
The proposition follows from \eqref{hung} and \eqref{uju}. 
\end{proof}

\begin{theorem}
    \label{exp}
    Let $u_0\in H^s(\R)\cap L^2(e^{2bx}dx)$ for $s\ge 0$ and $b>0$. There exist $T>0$ and a unique solution $u\in C([0,T];H^s(\R)\cap L^2(e^{2bx}dx))$ to the IVP \eqref{kdvks} satisfying \eqref{11} with $T$ instead of $T_0$. Moreover, for any $T'\in (0,T)$ there is a neighborhood $V$ of $u_0$ in $H^s(\R)\cap L^2(e^{2bx}dx)$ such that the data-solution map $u_0\mapsto u$ from $V$ into the class defined by \eqref{11} with $T$ instead of $T_0$ is Lipschitz.
\end{theorem}
\begin{proof}
    Fot $T>0$ consider $\lambda_8^T(w):=\|e^{bx}w\|_{L^\infty_TL^2_x}$. For $z>0$ to be defined, set \begin{equation*}
        W_z^T:=\{u\in C([0,T];H^s(\R)\cap L^2(e^{2bx}dx)) \mid \Theta^T(u)\le z\}, 
    \end{equation*}
where $\Theta^T(u):=\Lambda^T(u)+\lambda_6^T(u)+\lambda_8^T(u)$ with $\lambda_6^T(u)=A_6(T)\|\partial_x u\|_{L^2_TL^\infty_x}$ as in the proof of Theorem \ref{wLWP}.
Let us see that the integral operator $\Psi_{u_0}: W^T_z\to W^T_z$ defines a contraction in $W^T_z$ for some $T,z>0$.

Note for $0<t\le T$ that 
\begin{equation*}
    c(\eta,b)e^{-t\delta}(1+e^t)\le c(\eta,b)(1+e^T)^2:=c_*,
\end{equation*}
where $\delta$ is as in Proposition \ref{trocae}. Moreover, note $c_*$ can be taken independent on $T$ when $T\le1$.  From Proposition \ref{trocae} we get for $t\in[0,T]$ that
\begin{equation}\label{am1}
    \|e^{bx}V(t)u_0\|_{L^2_x}\le c_*\|e^{bx}u_0\|_{L^2_x}.
\end{equation}
Similarly, 
\begin{equation}\label{am2}
    \begin{split}
        \int_0^t\|e^{bx}V(t-t')u\partial_x u(t')\|_{L^2_x}dt'&\le c_*\int_0^t\|e^{bx}u\partial_x u(t')\|_{L^2_x}dt'\\&\le c_*T^{1/2}\|e^{bx}u\|_{L^\infty_TL^2_x}\|\partial_x u\|_{L^2_TL^\infty_x}\\&\le c_*T^{1/2}A_6\lambda_6^T(u)\lambda_8^T(u).
    \end{split}
\end{equation}
Thus, from \eqref{am1} and \eqref{am2} we have
\begin{equation}
    \lambda_8^T(\Psi(u))\le c_*\left\{\|e^{bx}u_0\|_{L^2_x}+T^{1/2}A_6\left(\Theta^{T}(u)\right)^2 \right\}.\label{e1}
\end{equation}

Also, we already know from \eqref{W1} that 
\begin{equation*}
    \lambda_6^T(u)\le c(\eta)\left\{\|u_0\|_{L^2_x}+A_2\Lambda^T(u)^2 \right\}.
\end{equation*}

Therefore
\begin{equation*}
    \Theta^T(\Psi(u))\le c(\eta,b,T)\left\{\|u_0\|_{H^s}+\|e^{bx}u_0\|_{L^2_x}+(A_2+T^{1/2}A_6)\left(\Theta^{T}(u)\right)^2  \right\}.
\end{equation*}

Set $z:=2c(\eta,b)\left\{\|u_0\|_{H^s}+\|e^{bx}u_0\|_{L^2_x}\right\}$. There exists $\tilde{T}>0$ such that $$(A_2+\tilde{T}^{1/2}A_6)(\Tilde{T})<\frac{1}{2c(\eta,b)z}.$$
Take $T:=\min\{1,\tilde{T}\}$. We conclude $\Theta^T(\Psi(u))\le z.$
The rest of the proof follows from standard-by-now arguments. 

\end{proof}

\section*{Acknowledgments}
This work was supported in part by CNPq, Conselho Nacional de Desenvolvimento Científico e Tecnológico - Brazil and INCTMat, Instituto Nacional de Ciência e Tecnologia de Matemática-Brazil.


\end{document}